\documentclass[12pt]{article}
\usepackage{amsmath,amsfonts,amssymb,amsthm}

\def\o{\omega}
\def\s{\sigma}

\def\CC{{\mathbb C}}

\def\ZZ{{\mathbb Z}}

\def\la{\langle}
\def\ra{\rangle}
\def\<{\langle}
\def\>{\rangle}

\begin{document}

\newcommand{\aut}{\mathrm{Aut}}
\newcommand{\M}{\mathbb{M}}
\newcommand{\dih}[2]{DIH_{#1}(#2)}
\newcommand{\drt}[1]{\frac{1}{\sqrt{2}} #1}
\newcommand{\drtp}[1]{  \frac{1}{\sqrt{2}}\left( #1\right)}
\newcommand{\drtpp}[1]{ \left( \frac{1}{\sqrt{2}}\left( #1\right)\right)}
\newcommand{\mspan}{\mathrm{span}}
\newcommand{\Rspan}[1]{\mathrm{span}_{\mathbb{R}}(#1)}
\newcommand{\Dih}[1]{DIH_{#1}}

\newtheorem{thm}{Theorem}[section]
\newtheorem{prop}[thm]{Proposition}
\newtheorem{lem}[thm]{Lemma}
\newtheorem{rem}[thm]{Remark}
\newtheorem{coro}[thm]{Corollary}
\newtheorem{conj}[thm]{Conjecture}
\newtheorem{de}[thm]{Definition}
\newtheorem{hyp}[thm]{Hypothesis}

\newtheorem{nota}[thm]{Notation}
\newtheorem{ex}[thm]{Example}
\newtheorem{proc}[thm]{Procedure}  

\newtheorem{tp}[thm]{Topic}  

\begin{center}\end{center} 

\begin{center}
{\Large  \bf  Research topics in finite groups and vertex algebras} 

\bigskip
20 March, 2019

\bigskip

Robert L. Griess Jr.

Department of Mathematics,

University of
Michigan,

Ann Arbor, MI 48109-1043  USA

{\tt rlg@umich.edu}

\bigskip

{\it Dedicated to Geoffrey Mason}
\bigskip

\end{center}

\tableofcontents

\section{Introduction}

In this article, we suggest a few projects for studying VOAs with emphasis on finite group viewpoints.   The world of VOAs is quite large and mostly unexplored.   It may be worthwhile to cultivate connections between the theories of finite groups and VOAs since many interesting VOAs are associated to finite groups. 

Some familiar VOAs, such as lattice type, have been extensively studied, both by Lie theoretic and discrete methods.   In contrast, the Monster VOA  has no nonzero derivations \cite{FLM} and much of its study involves properties of its finite automorphism group, the Monster.    

We will list a series of topics which seem worthwhile for research.  In a sense, this article is a successor to \cite{G-gnavoa1}.  

\subsection{A few definitions}
For basic background on vertex operators, see \cite {FLM}.   An {\it Ising vector} in a VOA  is a conformal vector of central charge $\frac 12$ which generates  a simple subVOA; it defines a {\it Miyamoto involution} (MI) on the VOA \cite{MI}.   A VOA has {\it CFT type} if it has the form $\bigoplus_{n\ge 0}V_n$ so that $dim(V_0)=1$.   A VOA has {\it OZ type} (or is an {\it OZVOA}) if it has CFT type and $V_1=0$ \cite{G-gnavoa1} (``OZ'' stands for one-dimensional, zero-dimensional).   Integral forms (IF) for a VOA are defined in \cite{DG-if}.     

Group theoretic notation will be generally consistent with \cite{Gor,GRqs}.

\section{Derivations}

For background, see \cite{DGfgvoa}.   

\begin{tp}
Find general conditions on a  VOA which make every derivation inner.    
\end{tp}

Such a result would be an analogue of the property that every derivation is inner for finite dimensional complex semisimple Lie algebras.   

\begin{tp}
Find general conditions on a VOA which make a VOA have non-inner derivations.   
Furthermore, find general upper and lower bounds on the dimension of the space of derivations in terms of dimensions of finite dimensional homogeneous subspaces which generate.   
\end{tp}

Results of this type might be analogous to results for finite groups about the size of the outer automorphism group.   
The latter issue can  involve first cohomology groups \cite{KG}.   
See also the G\" aschutz theorem that 
most finite $p$-groups have nontrivial outer automorphism group \cite{WG} and p.403 \cite{Hu1}.

\section{The $\o$-transposition property}  

\begin{tp} Find analogues in VOA theory of classifications in finite group theory by special generators, meaning the subgroups generated by two such elements are restricted.   See \cite{F0,F1,Timm,ThQP,Ho1,Ho2,Ho3,Ho4,Ho5,Ho6,Ho7,Asch-odd,Asch-odd-2,Asch-Hall}   
\end{tp}

The Miyamoto Involutions \cite{MI} are among the most accessible automorphisms of VOAs.    They do not exist for every VOA but they do exist in many cases where the VOA has finite automorphism groups.   In case the MIs have $\s$-type, they satisfy the famous 3-transposition condition of Fischer.   In case the MIs have $\tau$-type, and the VOA has a real form with positive definite natural bilinear form \cite{Sa}, the MIs satisfy a 6-transposition property (this is the case for $2A$ involution of the Monster action on the Monster VOA).   

\begin{tp}
Define other general types of involutions on families  of VOAs and prove that they satisfy some limited set of orders of products of pairs.   
\end{tp}

See \cite{HLY,HLY2,LY} for related results. 

\section{Automorphism groups}

\begin{tp}   Given a finite group $G$, is there a VOA whose automorphism group is $G$, or essentially $G$?    Better: is there such a VOA of CFT type?      
\end{tp}

The meaning of ``essentially'' here is $G$ occurs as $H/K$, where $H, K$ are normal subgroups of $Aut(V)$, for some VOA $V$,  and where $K$ and $Aut(V)/H$ are solvable and are small relative to $G$.   The first example of this is probably the case where $G$ is the Monster; the Moonshine VOA  has $G$ as its automorphism group \cite{FLM}.    There are recent results for $G\cong Sym_n$ \cite{LSY}.

If $V$ is a VOA of CFT and $V_1 \ne 0$, the (infinite) complex Lie group associated to the Lie algebra $V_1$ acts on $V$ as automorphisms, though not always faithfully.   An example is the Heisenberg VOA $M(1)$ of central charge 1, whose automorphism group is cyclic of order 2.   Rationality may be relevant \cite{DGfgvoa}.  

\begin{tp}   
Same as the previous topic except replace VOA over the complex numbers by a VA over a finite field.
\end{tp}

The answer is positive for the series of Chevalley-Steinberg groups defined over a given finite field \cite{GL-chst}.    The answer is unknown for almost all of the remaining families for groups of Lie type, the Suzuki groups $Sz(q)={}^2B_2(q)$ for $q\ge 8$ an odd power of $2$ and for the Ree groups
${}^2F_4(q)'$ for $q$ an odd power of 2 (note commutator subgroup at $q=2$) and ${}^2G_2(q)$ for $q\ge 27$ an odd power of 3 (${}^2G_2(3)\cong SL(2,8):3$).

\begin{tp} 
The twisted lattice type VOA for Barnes-Wall lattice of rank 32 has automorphism group of the form $2^{27}.E_6(2)$ \cite{Sh}.        
This group has the form of a parabolic subgroup of $E_7(2)$.    Find VOAs whose automorphism groups are essentially $E_6(2), E_7(2), E_8(2)$ and some of their parabolic subgroups.   
\end{tp} 

    If $V$ is the VOA of \cite{Sh}, $Aut(V)\cong 2^{27}{:}E_6(2)$.    Let  $A:=O_2(Aut(V))\cong 2^{27}$, then we have an action of $E_6(2)$ as automorphisms of $V^A$.   Possibly, $Aut(V^A)\cong E_6(2)$ or $E_6(2){:}2$.       
See \cite{G156} and Example 3.2 of \cite{G-gnavoa1}.

\section{Baer-Suzuki theorem for MIs}  

The Baer-Suzuki theorem says that in a finite group, $G$, if $p$ is a prime number and $x \in G$ has order a power of $p$, then $x \in O_p(G)$ if and only if $\la x, y \ra$ is a $p$-subgroup for all conjugates $y$ of $x$.   See \cite{Gor,AL}.   

We generalize a bit with this hypothesis: $p$ is a prime number, $K, L$ are conjugacy classes in the finite group $G$ with the property that $\la x, y \ra$ is a $p$-group, for all $x\in K, y\in L$.

The Baer-Suzuki theorem says that if $K=L$, $K$ is contained in $O_p(G)$.  If $K\ne L$,  $K$ may not be contained in $O_p(G)$.

For example, in $Sym_6$, take $p=2$ we let $K$ be the conjugacy class of $(12)$ and let $L$ be the conjugacy class of $(12)(34)(56)$. Then $\la x, y\ra$ is a 2-group for all $x\in K, y\in L$.    There may be counterexamples for $p=3$ in $G_2(3)$, but we do not know a reference.     

\begin{tp}  Is there a proof by VOA theory that the conclusion of the  Baer-Suzuki theorem holds for a conjugacy class of Miyamoto involutions?       
\end{tp}

\begin{tp}
   Is there an analogue of the Baer-Suzuki theorem which might be applicable to (short) words in MIs?    
\end{tp}

\section{Real, complex and quaternion reflection groups}

Finite real and complex reflection groups have been classified and studied a long time ago.  There is a less well-known classification of quaternionic reflection groups \cite{CQ}, whose list of conclusions includes double covers of alternating groups and the double cover of the Hall-Janko group.    

\begin{tp} Is there a class of VOAs which are quaternionic vector spaces and in which there a class of automorphisms of VOAs which behave like quaternionic reflections?  Are there such VOAs which  realize the finite groups generated by quaternionic reflections groups as their automorphism groups.  
\end{tp}

\section{Fixed point subVOAs}

\begin{tp}
Start with $V$ a familiar VOA, $F$ a finite subgroup of $Aut(V)$ and study $V^F$.   
Determine  properties of $V^F$, such as automorphism group, $C_2$, regularity, etc.   Decide whether  $V^F$ is isomorphic to a familiar VOA.  
One may continue by taking a finite subgroup $S$ of $Aut(V^F)$, then study $(V^F)^H$, etc.   
 \end{tp}

Parts of the above program were done for $V=V_{A_1}$ \cite{DGrank1ad,DGrank1e}.  See also the recent \cite{CM,HLY, HLY2,LY}.
The cases $V=V_L$ where $L$ is a root lattice of type ADE are likely to be interesting.   
There are many cases of finite group $F$ in $Aut(V)$ to try.   See the survey \cite{GRqs} and \cite{GRqs!}.   Note that there are examples in \cite{GRqs} for which $V_1^F=0$.   This means that $V^F$ is an OZVOA and in particular $V_2^F$ is a commutative algebra.

\subsection{An example}

Let $V=V_{E_8}$ and let $F$ be the {\it Borovik group}, i.e. a subgroup $B$ of $Aut(V)\cong E_8(\CC )$ which has these properties: 

(i) $B$ contains a normal subgroup $B_5 \times B_6$, where $B_k\cong Alt_k$ and $C_B(B_5)=B_6$ and $C_B(B_6)=B_5$;      
  
(ii) By conjugation on $B_5 \times B_6$, $B$ embeds as an index 2 subgroup of $Aut(B_5)\times Aut(B_6)$;

(iii) $B$ contains a subgroup isomorphic to $Alt_5 \times Sym_6$; and 

(iv) $B$ is not a nontrivial direct product.   

\bigskip

This is a Lie primitive subgroup of $E_8(\CC )$ (i.e., not contained in a positive dimensional Zariski-closed subgroup except for $E_8(\CC )$).   See \cite{Boro,CG} for characterization of $B$ and \cite{FrGr} for the $E_8(\CC )$-conjugacy classes represented in $B$.    

In particular, $V_1^B=0$, $dim(V_2^{B_5})=16$ and $dim(V_2^{B_6})=10$ \cite{G-gnavoa1} .   Since $B/B_5\cong Aut(Alt_6)$ acts as automorphisms on the commutative algebra $V_2^{B_5}$, a reasonable guess would be that $V_2^{B_5}\cong Mat_{4\times 4}^+$, a Jordan algebra.   Since $B/B_6\cong Sym_5$ acts as automorphisms on the commutative algebra $V_2^{B_6}$, a reasonable guess would be that $V_2^{B_5}\cong SymMat_{4\times 4}$, a Jordan algebra.   

\begin{tp} Identify the fixed point subalgebras $V_2^{B_5}$ and $V_2^{B_6}$.   Determine their automorphism groups.   
\end{tp}

\begin{tp}   For $F=B_5$ or $B_6$, let $V_2^F$ be the fixed point subalgebra in $V_2$.   Find Ising vectors (conformal vectors of central charge 1/2) in $V_2$ whose Miyamoto involutions give the $2A$-involutions of $B$; see \cite{FrGr}.    
\end{tp}

We remark that the function which sends an Ising vector to its Miyamoto involution is in general not a monomorphism \cite{GH}.  Given a $2A$-involution of $B$, it is possible that some Ising vectors which represent it are in $V_2^F$, or that none are.

\section{Integral form classifications for $V$ and $V_2$}  

Studies of integral forms in VOAs are not (yet) numerous.   See \cite{GIF} for a survey.    Most results are existence proofs.    

\begin{tp}  Give existence proofs of integral forms for more VOAs of interest, and give uniqueness proofs (say for maximal integral forms).  
\end{tp}

When the automorphism group of a VOA is a positive dimensional algebraic group, an integral form has an infinite orbit under the automorphism group, so can not be unique.  
A uniqueness program must be sensitive to the fact that different integral forms may have different signatures (as real quadratic forms).    

\begin{tp}  Determine a suitable notion of equivalence for group-invariant integral forms (or maximal integral forms) in a given VOA, and give conditions under which integral forms have finitely many equivalence classes (and even one equivalence class).
\end{tp}

A step in this direction was obtained in the thesis of Simon \cite{Simon}, who classified maximal integral forms in $S$ which are invariant under certain finite groups.   Here, $S$ is one of the nine dihedral algebras determined by Sakuma \cite{Sa}.  These algebras are naturally labeled by nodes of the extended $E_8$-diagrams.             This $S$ occurs as $V_2$, where $V$ is a VOA generated by a pair of Ising vectors corresponding to one of the  nine configurations in \cite{Sa}.    

\begin{tp}  Determine the maximal, group-invariant integral forms in a given dihedral VOA, $V$.   Is it true that such an integral form is generated by one of the integral forms in $V_2$ as classified in \cite{Simon}?
\end{tp}

\section{Uniqueness of dihedral VOA}

Suppose that the VOA $V$ has OZ type.

\begin{tp}  To what extent is $V$  is determined by its degree 2 summand, $V_2$?    
\end{tp}

A necessary condition for $V$ to be uniquely determined by $V_2$ is that it be generated by its degree 2 summand.   

\begin{tp}   When $V_2$ is one of the nine Sakuma dihedral algebras and $V_2$ generates $V$, is $V$ unique?
\end{tp}

The answer is yes for at lease one case of the Sakuma dihedral algebras. See \cite{DJY}.   

The Sakuma classification \cite{Sa} follows from a more general later result \cite{HRS}.  

\section{Integral Forms}

Existence proofs for integral forms in VOAs have begun to appear in recent literature \cite{GIF}.   Classifications of integral forms are fewer.   The thesis of Simon \cite{Simon} classifies certain maximal group-invariant integral forms in the nine degree 2 Sakuma dihedral algebras, and in a few other commutative nonassociative algebras.   

\begin{tp} Extend the above classifications of IFs in Sakuma algebras to IFs in full dihedral VOAs.
\end{tp}

\begin{tp}  Suppose that the $V$ is holomorphic.   Is there self-dual integral form in $V$ unique, up to equivalence.    In particular, does the Moonshine VOA have a unique self-dual integral form?  
\end{tp}

In \cite{Car}, Carnahan proposes four integral self-dual forms for the Moonshine VOA and conjectures that they ought to be equal.

\section{Module theory for Chevalley-Steinberg groups }

Let $\mathbb F \le \mathbb E$ be a finite degree field extension and $\Phi$ an indecomposable root system.   In \cite{GL-chst},   Griess and Lam gave constructions for classical VAs of type $\Phi$ over $\mathbb E$  and of the related Chevalley-Steinberg groups 
associated to $\mathbb F \le \mathbb E$.   These groups are essentially the full automorphism groups for these VAs.   
 
For the positive characteristic case, the finite dimensional modules for these groups are generally not  completely reducible.      

\begin{tp} For each graded piece of the VA, determine the composition factors with respect to the associated Chevalley-Steinberg group and determine a direct sum decomposition into indecomposable modules.  
\end{tp}

\begin{tp} By examining the graded pieces of the VA, deduce lower bounds for dimensions of $Ext^n (A,B)$ for nonzero modules $A, B$ for the automorphism group and for its socle.    
 \end{tp}

Constructing nonsplit extensions of one module by another can be technically difficult.      
The point here is that the VA graded space gives at once a vast collection of modules to study.   

\section{Acknowledgments}

Our talk given at the Mason fest 2018, a survey of integral form results in vertex operator algebras,  is represented in articles already committed to other publications \cite{GIF,DGdet-if}.  

Work on the present article was supported by funds from University of Michigan.  

We thank the referee for comments and references.   
 


\begin{thebibliography}{GRC99}   


\bibitem{AL} 
MR0289622  Alperin, J.; Lyons, R. On conjugacy classes of $p$-elements. J. Algebra 19 1971 536-537. (Reviewer: D. K. Friesen)

\bibitem{Asch-odd} MR0310058  Aschbacher, Michael On finite groups generated by odd transpositions. I. Math. Z. 127 (1972), 45-56. (Reviewer: R. W. Carter) 20D99

\bibitem{Asch-odd-2} 
MR0327900  Aschbacher, Michael On finite groups generated by odd transpositions. II, III, IV. J. Algebra 26 (1973), 451-459; ibid. 26 (1973), 460-478; ibid. 26 (1973), 479?491. (Reviewer: R. W. Carter) 20D99

\bibitem{Asch-Hall}  MR0311765  
Aschbacher, Michael; Hall, Marshall, Jr.
Groups generated by a class of elements of order 3.
J. Algebra 24 (1973), 591-612.
20D99

\bibitem{Boro} MR1066315 (91h:20071) 
Borovik, A. V.
The structure of finite subgroups of simple algebraic groups. (Russian)
Algebra i Logika 28 (1989), no. 3, 249--279, 366; translation in
Algebra and Logic 28 (1989), no. 3, 163
 
\bibitem{Car} Carnahan, Scott: 
Four self-dual integral forms of the moonshine module; 
arXiv:1710.00737  

\bibitem{CM} 
 Scott Carnahan, Masahiko Miyamoto, Regularity of fixed-point vertex operator subalgebras, arXiv:1603.05645.

\bibitem{CG} MR1200322  Cohen, Arjeh M.; Griess, Robert L., Jr. Nonlocal Lie primitive subgroups of Lie groups. Canad. J. Math. 45 (1993), no. 1, 88?103. (Reviewer: Martin W. Liebeck) 22E10

\bibitem{CQ} 579063 Cohen, Arjeh M. Finite quaternionic reflection groups. J. Algebra 64 (1980), no. 2, 293--324. (Reviewer: J. D. Dixon) 20H15 

\bibitem{DGrank1ad} MR1644007   Dong, Chongying; Griess, Robert L., Jr. Rank one lattice type vertex operator algebras and their automorphism groups. J. Algebra 208 (1998), no. 1, 262?275. (Reviewer: Hai Sheng Li) 17B69



\bibitem{DGrank1e}  MR1700522  Dong, Chongying; Griess, Robert L., Jr.; Ryba, Alex Rank one lattice type vertex operator algebras and their automorphism groups. II. E-series. J. Algebra 217 (1999), no. 2, 701--710. (Reviewer: Hai Sheng Li) 17B69

\bibitem{DGfgvoa} MR1914063  Dong, Chongying; Griess, Robert L., Jr. Automorphism groups and derivation algebras of finitely generated vertex operator algebras. Michigan Math. J. 50 (2002), no. 2, 227--239. (Reviewer: Pavel S. Kolesnikov) 17B69

\bibitem{DG-if} MR2928458  Dong, Chongying; Griess, Robert L., Jr. Integral forms in vertex operator algebras which are invariant under finite groups. J. Algebra 365 (2012), 184--198. (Reviewer: Drazen Adamovic) 17B69 (17B68)



\bibitem{DGdet-if} Chongying Dong and Robert L. Griess, Jr.;  
Determinants for integral forms in lattice type vertex operator algebras, submitted; arxiv:1903.08210.

\bibitem{DJY}   Chongying Dong; Xiangyu Jiao;
Nina Yu; 6A-Algebra and its representations; preprint.    


\bibitem{F0}  Fischer, Bernd; University of Warwick Notes, about 1968

\bibitem{F1}  MR0294487  Fischer, Bernd;  Finite groups generated by 3-transpositions. I. Invent. Math. 13 (1971), 232--246. (Reviewer: A. R. Camina) 20D99

\bibitem{FLM}  Frenkel, Igor; Lepowsky, James; Meurman, Arne;  Vertex operator algebras and the Monster. Pure and Applied Mathematics, 134. Academic Press, Inc., Boston, MA, 1988. liv+508 pp. ISBN: 0-12-267065-5 (Reviewer: Kailash C. Misra) 17B65 (17B67 20D08 81D15 81E40)

\bibitem{FrGr}
MR1620658 (2000a:20035) 
Frey, Darrin D.; Griess, Robert L., Jr.(1-MI)
Conjugacy classes of elements in the Borovik group. (English summary)
J. Algebra 203 (1998), no. 1, 226--243.
20D06
  
 
\bibitem{WG}  MR0193144  Gasch\" utz, Wolfgang:  Nichtabelsche $p$-Gruppen besitzen \"aussere $p$-Automorphismen. (German) J. Algebra 4 1966 1?2. (Reviewer: I. Reiner) 20.40

\bibitem{Gor}  MR0231903  Gorenstein, Daniel :  Finite groups. Harper \& Row, Publishers, New York-London 1968 xv+527 pp. (Reviewer: W. Feit) 20.25;
MR0569209  Gorenstein, Daniel Finite groups. Second edition. Chelsea Publishing Co., New York, 1980. xvii+519 pp. ISBN: 0-8284-0301-5 20-02 (20Dxx)


\bibitem{G-gnavoa1} MR2029791  Griess, Robert L., Jr. GNAVOA. I. Studies in groups, nonassociative algebras and vertex operator algebras. Vertex operator algebras in mathematics and physics (Toronto, ON, 2000), 71--88, Fields Inst. Commun., 39, Amer. Math. Soc., Providence, RI, 2003. (Reviewer: Hai Sheng Li) 17B69 (20B25);

MR1980921 Indexed
Griess, Robert L., Jr.(1-MI)
Summary of GNAVOA. I. Studies in groups, nonassociative algebras and vertex operator algebras. (English summary)
Algebraic combinatorics (Japanese) (Kyoto, 2001).
Surikaisekikenkyu sho Kokyuroku No. 1299 (2003), 1--6.
17B69




\bibitem{GIF} Robert L. Griess, Jr.; 
Integral forms 
in vertex operator algebras, a survey
{\it  Lecture at Ischia Group Theory Meeting, 
20 March, 2018}; submitted; arxiv:1903.08292  

\bibitem{GH} 
MR1998606  Griess, Robert L.; H\" ohn, Gerald Virasoro frames and their stabilizers for the $E_8$ lattice type vertex operator algebra. J. Reine Angew. Math. 561 (2003), 1--37. (Reviewer: Hai Sheng Li) 17B69 (17B67)

\bibitem{GL-chst}   MR3456026  Griess, Robert L., Jr.; Lam, Ching Hung Groups of Lie type, vertex algebras, and modular moonshine. Int. Math. Res. Not. IMRN 2015, no. 21, 10716?10755. (Reviewer: Robert Harold McRae) 17B69 (17B20 17B50 20D06 20D08);
MR3356595  Griess, Robert L., Jr.; Lam, Ching Hung Groups of Lie type, vertex algebras, and modular moonshine. Electron. Res. Announc. Math. Sci. 21 (2014), 167--176. 20D05 (17B69)

\bibitem{GRqs} MR1879514  Griess, Robert L., Jr.; Ryba, A. J. E. Classification of finite quasisimple groups which embed in exceptional algebraic groups. J. Group Theory 5 (2002), no. 1, 1--39. (Reviewer: Martin W. Liebeck) 20G20 (20D06 20E28);

\bibitem{GRqs!} 
MR1653177  Griess, Robert L., Jr.; Ryba, A. J. E. Finite simple groups which projectively embed in an exceptional Lie group are classified! Bull. Amer. Math. Soc. (N.S.) 36 (1999), no. 1, 75--93. (Reviewer: D. Wales) 20G20 (20C15 20E32 20G05)

\bibitem{G156} MR1650633  Griess, Robert L., Jr. A vertex operator algebra related to $E_8$ with automorphism group $O^+(10,2)$. The Monster and Lie algebras (Columbus, OH, 1996), 43--58, Ohio State Univ. Math. Res. Inst. Publ., 7, de Gruyter, Berlin, 1998. (Reviewer: Koichiro Harada) 17B69 (20D08)

\bibitem{KG}   Karl Gruenberg, 
Cohomological Topics in Group Theory (Springer Lecture Notes in Mathematics; 143), 1970.    

\bibitem{HRS} 
Hall, J. I.(1-MIS); Rehren, F.(4-BIRM-M); Shpectorov, S.(4-BIRM-M): Universal Axial Algebras and a Theorem of Sakuma, about 31 pages; to appear in Journal of Algebra.  
17B69


\bibitem{HLY} H\"ohn, Gerald; Lam, Ching Hung; Yamauchi, Hiroshi McKay?s E7 observation on the Baby Monster. Int. Math. Res. Not. IMRN 2012, no. 1, 166?212.

\bibitem{HLY2} H\"ohn, Gerald; Lam, Ching Hung; Yamauchi, Hiroshi McKay?s E6 observation on the largest Fischer group. Comm. Math. Phys. 310 (2012), no. 2, 329?365.


\bibitem{Hu1}   MR0224703 Huppert, B. Endliche Gruppen. I. (German)
Die Grundlehren der Mathematischen Wissenschaften, Band 134 Springer-Verlag, Berlin-New York 1967 xii+793 pp.


\bibitem{Ho1}
MR2611698 Thesis Ho, Chat-Yin ON THE CLASSICAL CASE OF THE QUADRATIC-PAIRS FOR 3 WHOSE ROOT-GROUP HAS ORDER 3. Thesis (Ph.D.) The University of Chicago. 1972. (no paging), ProQuest LLC


\bibitem{Ho2} 
MR0344341  Ho, C. Y. On the quadratic pair whose root group has order 3. Bull. Inst. Math. Acad. Sinica 1 (1973), no. 2, 155--180. (Reviewer: Bhama Srinivasan) 20D99


\bibitem{Ho3}
MR0384923  Ho, Chat Yin Quadratic pairs for 3 whose root group has order greater than 3. I. Comm. Algebra 3 (1975), no. 11, 961--1029. (Reviewer: Stephen D. Smith) 20D05


\bibitem{Ho4} 
MR0414737  Ho, C. Y. Chevalley groups of odd characteristic as quadratic pairs. J. Algebra 41 (1976), no. 1, 202--211. (Reviewer: Bruce Cooperstein) 20H20



\bibitem{Ho5}
MR0417273  Ho, C. Y. Quadratic pairs for odd primes. Bull. Amer. Math. Soc. 82 (1976), no. 6, 941--943. (Reviewer: Peter M. Neumann) 20D05

\bibitem{Ho6}
MR0422404  Ho, Chat-Yin On the quadratic pairs. J. Algebra 43 (1976), no. 1, 338--358. (Reviewer: Stephen D. Smith) 20D05 (20G40)

\bibitem{Ho7}
MR0778986  Baumann, B.; Ho, C. : Linear groups generated by a pair of quadratic action subgroups. Arch. Math. (Basel) 44 (1985), no. 1, 15--19. 20D05 (20D20) 


\bibitem{LSY} Lam, Ching Hung; Sakuma, Shinya; Yamauchi, Hiroshi Ising vectors and automorphism groups of commutant subalgebras related to root systems. Math. Z. 255 (2007), no. 3, 597 ?626.

\bibitem{LY} Lam, Ching Hung; Yamauchi, Hiroshi On 3-transposition groups generated by ?- involutions associated to c = 4/5 Virasoro vectors. J. Algebra 416 (2014), 84?121.
 


\bibitem{MI} MR1367861  Miyamoto, Masahiko:  Griess algebras and conformal vectors in vertex operator algebras. J. Algebra 179 (1996), no. 2, 523--548. (Reviewer: Chong Ying Dong) 17B69 (17B68)


\bibitem{M1} Miyamoto, Masahiko, $C_2$-cofiniteness of cyclic-orbifold models. Comm. Math. Phys. 335 (2015), no. 3, 1279?1286.

\bibitem{M2} Miyamoto, Masahiko, $C_2$-cofiniteness of orbifold models for finite groups, arXiv:1812.00570.




\bibitem{Sa} MR2347298  Sakuma, Shinya:   6-transposition property of $\tau$ -involutions of vertex operator algebras. Int. Math. Res. Not. IMRN 2007, no. 9, Art. ID rnm 030, 19 pp. (Reviewer: Ozren Per¨e) 17B69

\bibitem{Sh} MR3177874  Shimakura, Hiroki The automorphism group of the $\ZZ_2$-orbifold of the Barnes-Wall lattice vertex operator algebra of central charge 32. Math. Proc. Cambridge Philos. Soc. 156 (2014), no. 2, 343--361. (Reviewer: Alonso Castillo-Ramirez) 17B69

\bibitem{Simon} MR3641120  Simon, Gregory G. Automorphism-invariant Integral Forms in Griess Algebras. Thesis (Ph.D.) University of Michigan. 2016. 137 pp. ISBN: 978-1369-58993-1, ProQuest LLC

\bibitem{ThQP} MR0430043  
Thompson, J. G.
Quadratic pairs. Actes du Congr\` es International des Math\' ematiciens (Nice, 1970), Tome 1, pp. 375--376. Gauthier-Villars, Paris, 1971.
20C05 (22E99)

\bibitem{Timm} MR1852057 (2002f:20070) 
Timmesfeld, Franz Georg(D-GSSN)
Abstract root subgroups and simple groups of Lie type. (English summary)
Monographs in Mathematics, 95. Birkhäuser Verlag, Basel, 2001. xiv+389 pp. ISBN: 3-7643-6532-3
20G15 (20E42)


\end{thebibliography}
\end{document}